\documentclass[11pt]{article}
\usepackage{amsmath, amsfonts, amssymb}
\usepackage{graphicx, amsthm}
\usepackage{cite}
\usepackage{color}
\usepackage{amsmath}
\usepackage{enumitem}

\theoremstyle{definition}
\newtheorem{thm}{Theorem}[section]
\newtheorem{lem}[thm]{Lemma}

\newtheorem{cor}[thm]{Corollary}
\newtheorem{prop}[thm]{Proposition}

\newtheorem{defi}[thm]{Definition}
\theoremstyle{remark}
\newtheorem{rem}[thm]{Remark}
\numberwithin{equation}{section}
\theoremstyle{plain}

\newcommand{\R}{\mathbb{R}}

\newcommand{\ra}{\sqrt}
\newcommand{\ep}{\varepsilon}

\newcommand{\la}{\langle}
\renewcommand{\ra}{\rangle}

\begin{document}

\begin{center}
{\Large A partial answer to the Demyanov-Ryabova conjecture}\footnote{Research
partially supported by the grants BASAL PFB-03 and ECOS/CONICYT--ECOS/Sud
C14E06. The first author has also been partially supported by the grants
FONDECYT 1171854 (Chile) and MTM2014-59179-C2-1-P (MINECO of Spain and ERDF of
EU)}

\vspace{0.6cm}

\textsc{Aris Daniilidis \& Colin Petitjean}
\end{center}

\bigskip

\noindent\textbf{Abstract.} In this work we are interested in the
Demyanov--Ryabova conjecture for a finite family of polytopes. The conjecture
asserts that after a finite number of iterations (successive dualizations),
either a 1-cycle or a 2-cycle eventually comes up. In this work we establish a
strong version of this conjecture under the assumption that the initial family
contains ``enough minimal polytopes" whose extreme points are ``well placed".

\bigskip


\vspace{0.6cm}

\textit{Secondary} 49J52, 90C49.

\section{Introduction}

We call $\emph{polytope}$ any convex compact subset of $\mathbb{R}^{N}$ with a
finite number of extreme points. Throughout this work we consider a finite
family $\Re=\{\Omega_{1},\ldots,\Omega_{\ell}\}$ of polytopes of
$\mathbb{R}^{N}$ together with an operation which transforms the initial
family $\Re$ to a dual family of polytopes that we denote $\mathcal{F}(\Re).$
(Motivation and origin of this operation will be given at the end of the introduction).

\smallskip

Let us now describe the operation $\mathcal{F}$: let $\mathrm{ext}(\Omega)$
stand for the set of extreme points of the polytope $\Omega$ and let $S$
denote the unit sphere of $\mathbb{R}^{N}.$ Then given a family $\Re$ as
before, for any direction $d\in S$ and polytope $\Omega_{i}\in\Re$
($i\in\left\{  1,\ldots,\ell\right\}  $) we consider the set of $d$-active
extreme points of $\Omega_{i}$
\begin{equation*}
E(\Omega_{i},d):=\{x\in\mathrm{ext}(\Omega_{i}):\,\langle x,d\rangle
\ =\max\langle\Omega_{i},d\rangle\}. \label{d-active}
\end{equation*}
We associate to $d\in S$ the polytope
\begin{equation}
\Omega(d):=\mathrm{conv}\left(  \bigcup\limits_{\Omega_{i}\in\Re}
\,E(\Omega_{i},d)\right)  , \label{Omega(d)}
\end{equation}
that is, the polytope obtained as convex hull of the set of all $d$-active
extreme points (when $\Omega_{i}$ is taken throughout $\Re$). Since the set of
extreme points of all polytopes of the family $\Re$
\begin{equation}
E_{\Re}=\bigcup\limits_{\Omega_{i}\in\Re}\,\mathrm{ext}(\Omega_{i}) \label{E}
\end{equation}
is finite, the family of polytopes
\begin{equation}
\mathcal{F}(\Re):=\{\Omega(d):d\in S\} \label{F(R)}
\end{equation}
is also finite, hence of the same nature as $\Re$. We call $\mathcal{F}(\Re)$
the dual family of $\Re.$

\smallskip

Now starting from a given family of polytopes $\Re_{0}$, we define
successively a sequence of families $\{\Re_{n}\}_{n}$ by applying repeatedly
this duality operation (transformation) $\mathcal{F}$, that is, setting
$\Re_{n+1}:=\mathcal{F}(\Re_{n}),$ for all $n\in\mathbb{N}$. Since the
transformation $\mathcal{F}$ cannot create new extreme points, the sequence
\[
E_{\Re_{n}}=\bigcup\limits_{\Omega\in\Re_{n}}\,\mathrm{ext}(\Omega
)\quad\text{(extreme points of polytopes in }\Re_{n}\text{)}\qquad
n\in\mathbb{N}
\]
is nested (decreasing) and eventually becomes stable, equal to a finite set
$E.$ By a standard combinatorial argument, we now deduce that for some
$k\geq1$ and $n_{0}\geq0$ we necessarily get $\Re_{n}=\Re_{n+k}$ (and
$E_{\Re_{n}}=E$), for all $n\geq n_{0}$. Therefore, a $k$-cycle $(\Re_{n_{0}
},\Re_{n_{0}+1},\cdots,\Re_{n_{0}+k-1})$ is always formed. We are now ready to
announce the conjecture of Demyanov and Ryabova:

\begin{itemize}
\item \textbf{Conjecture} (Demyanov--Ryabova, \cite{demrya}). Let $\Re_{0}$ be
a finite family of polytopes in $\mathbb{R}^{N}.$ Then for some $n_{0}
\in\mathbb{N}$ we shall have $\Re_{n_{0}}=\Re_{n_{0}+2}$.
\end{itemize}

In other words, after some threshold $n_{0}$ the sequence
\[
\Re_{0},\quad\Re_{1}=\mathcal{F}(\Re_{0}),\cdots,\quad\;\Re_{n+1}
=\mathcal{F}(\Re_{n}),\cdots
\]
stabilizes to either a $1$-cycle (self-dual family $\Re_{n}=\mathcal{F}
(\Re_{n})=\Re_{n+1}$) or to a $2$-cycle (reflexive family $\Re_{n}
=\mathcal{F}(\mathcal{F}(\Re_{n}))=\Re_{n+2}$) for $n\geq n_{0}$. In
\cite{demrya}, the authors carried out generic numerical experiments over two
hundred families of polytopes, where only $2$-cycles eventually arise. Notice
however that one can construct particular examples where a $1$-cycle is
formed. In all known cases, the initial family $\Re_{0}$ ends up, after finite
iterations, to a reflexive one.

\smallskip

Besides the recorded numerical evidence, there is still no proof of this
conjecture. The only known result in this direction is due to \cite{sang}. In
that work, the author establishes the conjecture under the additional
assumption that the set $E_{\Re_{0}}$ of extreme points of the initial family
$\Re_{0}$ is affinely independent.

\smallskip

Before we state and prove our main result, let us mention that in
$1$--dimension the conjecture is trivially true.

\begin{prop}
[The conjecture is true in $1$--dim]\label{Prop-dim-1} Let $\Re_{0}$ be a
finite family of closed bounded intervals of $\mathbb{R}$. Then $\Re_{1}
=\Re_{3}$.
\end{prop}

\noindent\textbf{Proof.} Let us denote $\{I_{1},\ldots,I_{\ell}\}$ the
elements of $\Re_{0}$ with $I_{j}=[a_{j},b_{j}],$ $j\in\{1,\ldots,\ell\}$.
Since the unit sphere $S=S_{\mathbb{R}}=\{1,-1\}$ consists of only two
directions, the construction of the dual family $\Re_{1}=\mathcal{F}(\Re_{0})$
is very simple. To this end, we set $a_{-}:=\min_{i\in\{1..\ell\}}a_{i}$,
$a_{+}:=\max_{i\in\{1..\ell\}}a_{i}$, $b_{-}:=\min_{i\in\{1..\ell\}}b_{i}$,
$b_{+}:=\max_{i\in\{1..\ell\}}b_{i}$. This leads to the family
\[
\Re_{1}=\{\Omega_{1}(-1),\Omega_{1}(1)\}=\{[a_{-},a_{+}],[b_{-},b_{+}]\}.
\]
The construction of $\Re_{2}=\mathcal{F}(\Re_{1})$ is even simpler, since we
only have two intervals (polytopes) to consider. We actually have
\[
\Re_{2}=\{\Omega_{2}(-1),\Omega_{2}(1)\}=\{[a_{-},b_{-}],[a_{+},b_{+}]\}.
\]
It now suffices to compute $\Re_{3}$ and obtain directly that $\Re_{1}=\Re
_{3}$. (Notice that if it happens $a_{+}=b_{-}$ then we actually get a
$1$-cycle: $\Re_{1}=\Re_{2}$.)\hfill$\square$

\bigskip

The extreme simplicity of the problem in dimension $1$ is due to the fact that
the family that arises after any new iteration has at most $2$ elements
(corresponding to the directions $1$ and $-1$ of the unit sphere
$S_{\mathbb{R}}$). The problem gets much more complicated though in higher
dimensions, where no prior efficient control on the cardinality of the
iterated families can be obtained (apart from an absolute combinatorial bound
on the number of all possible polytopes that can be obtained by convexifying
subsets of the prescribed set of extreme points $E$). We shall now treat this
general case.

\smallskip

Let $\Re_{0}$ be a finite family of polytopes in $\mathbb{R}^{N}$ ($N\geq 2$). We
denote by $E:=E_{\Re_{0}}\ $the set of extreme points of all polytopes of the
family, see (\ref{E}), by $R=\mathrm{card}(E)$ its cardinality and we set
\begin{equation*}
C:=\mathrm{conv}(E) \label{C}
\end{equation*}
its convex hull. Notice that every polytope $\Omega$ of the family $\Re_{0}$
(or of any family $\Re_{n}$ obtained after $n$-iterations, for every
$n\in\mathbb{N}$), is contained in $C.$ Let further
\[
r(\Omega):=\mathrm{card}(\Omega\cap E)
\]
denote the number of extreme points of the polytope $\Omega\in\Re_{0}$ and
set
\begin{equation}
r_{\min}:=\min_{\Omega\in\Re_{0}}r(\Omega). \label{r-min}
\end{equation}
We now state the main result of the paper.

\begin{thm}
[Main result]\label{hyp}Let $\Re_{0}$ be a finite family of polytopes in
$\mathbb{R}^{N}$ and $r_{\min}\in\{1,\ldots,R\}$ as in (\ref{r-min}). Then
$\Re_{1}=\Re_{3}$ (\textit{i.e.} a reflexive family occurs after one iteration) provided:

\begin{itemize}
\item[(H1)] $\forall\,x\in E,\,x\not \in \mathrm{conv}(E\backslash
\{x\})$ \quad(\textit{i.e.} each $x\in E$ is extreme in $C.$)

\item[(H2)] $\Re_{0}$ contains all $r_{\min}$-polytopes (that is, all
polytopes made up of $r_{\min}$ points of $E$).
\end{itemize}
\end{thm}

\smallskip

\begin{rem}
\textrm{(i)} Assumption (H$_{1}$) easily yields that the set of extreme points
remains stable from the very beginning, that is,
\begin{equation*}
E_{\Re_{n}}=E_{\Re_{0}}=E,\qquad\text{for all }n\in\mathbb{N}. 
\end{equation*}
Indeed pick $x\in E$ and $\mathrm{e}_{x}\in S$ which exposes $x$ in $C$. Let $\Omega
\in\Re_{0}$ be such that $x\in\mathrm{ext}(\Omega)$ (there is clearly at least
one such a polytope in $\Re_{0}$). Then $\mathrm{e}_{x}$ exposes $x$ in $\Omega$, that is
$x\in E(\Omega,\mathrm{e}_{x}).$ It follows readily that $x\in\Omega(\mathrm{e}_{x})\subset
E_{\Re_{1}}$ (see the definition in (\ref{Omega(d)})) and by a simple induction, $x\in E_{\Re_{n}}$, for every $n\geq
1$.\smallskip\newline\textrm{(ii) }Assumption (H$_{2}$) will be weakened in
the sequel.
\end{rem}

\textbf{Origin of the conjecture.} The initial motivation which eventually led to the
formulation of the above conjecture stems from the problem of stable
representation of positively homogeneous polyhedral functions as finite minima
of sublinear ones, or its geometric counterpart, the representation of a
closed polyhedral cone as finite union of closed convex polyhedral cones. Let
us recall that a function $f:\mathbb{R}^{N}\rightarrow\mathbb{R}$ is called
positively homogeneous provided $f(\lambda x)=\lambda f(x)$ for every
$x\in\mathbb{R}^{N}$ and $\lambda>0$. It is called sublinear (respectively,
superlinear) if it is positively homogeneous and convex (respectively,
concave). \smallskip

Following \cite{psh}, a sublinear function $\overline{g}:\mathbb{R}
^{N}\rightarrow\mathbb{R}$ is called an upper convex approximation of $f$ if
$\overline{g}$ majorizes $f$ on $\mathbb{R}^{N}$, that is, $\overline
{g}(x)\geq f(x),$ for every $x\in\mathbb{R}^{N}$. In the same way, a
superlinear function $\underline{g}$: $\mathbb{R}^{N}\rightarrow\mathbb{R}$ is
called a lower concave approximation of $f$ if $\underline{g}$ minorizes $f$
on $\mathbb{R}^{N}$, that is, $\underline{g}(x)\leq f(x)$ for all
$x\in\mathbb{R}^{N}$. Then we say that a set of sublinear functions
$E^*$ is an upper exhaustive family for $f$ if the following equality
holds for every $x\in\mathbb{R}^{N}$:
\begin{equation}
f(x)=\inf_{\overline{g}\in E^*}\overline{g}(x). \label{upper}
\end{equation}
Similarly, we say that a set of superlinear functions $E_*$ is a
lower exhaustive family for $f$ if the following equality holds for every
$x\in\mathbb{R}^{N}$:
\begin{equation}
f(x)=\sup_{\underline{g}\in E_*}\overline{g}(x). \label{lower}
\end{equation}
In \cite{rubinov} the authors established the existence of an upper exhaustive
family of upper convex approximations (respectively lower exhaustive family of
lower concave approximations) when $f$ is upper semicontinuous on
$\mathbb{R}^{N}$ (respectively lower semicontinuous). In particular, if $f$ is
continuous, the existence of both such families is guaranteed.

\smallskip

It is well known (see \cite{LHU,phelps} \emph{e.g.}) that a function
$\overline{g}$: $\mathbb{R}^{N}\rightarrow\mathbb{R}$ is sublinear if and only
if ${\overline{g}(x)=\max_{h\in\partial\overline{g}(0)}\langle x,h\rangle}$.
Using this fact we are able to restate (\ref{upper}) in the following way:
\[
f(x)=\inf_{\overline{g}\in E^*}\overline{g}(x)=\inf_{\overline{g}
\in E^*}\max_{h\in\partial\overline{g}(0)}\ \langle x,h\rangle
=\inf_{\bar{\Omega}\in\overline{\Re}}\max_{h\in\bar{\Omega}}\ \langle
x,h\rangle,
\]
where $\overline{\Re}=\{\partial\overline{g}(0):\overline{g}\in E^*\}$
is the family of subdifferentials of the sublinear functions $\overline{g}$
that represent $f$ and $\bar{\Omega}=\partial\overline{g}(0)$. In a similar
way, considering superlinear functions $\underline{g}$ (lower concave
approximations of $f$) and denoting by $\underline{\Re}=\{-\partial
(-\underline{g})(0):\underline{g}\in E_*\}$ the family of
superdifferentials $\underline{\Omega}=-\partial(-\underline{g})(0)$, we can
restate (\ref{lower}) as follows:
\[
f(x)=\sup_{\underline{\Omega}\in\underline{\Re}}\min_{h\in\underline{\Omega}
}\ \langle x,h\rangle.
\]
In case of a polyhedral function $f$ the exhaustive families $E^*$
and $E_*$ can be taken to be finite, with elements being polyhedral
functions ($\overline{g}$ and $\underline{g}$ respectively). In this case, the
corresponding families $\overline{\Re}$ and $\underline{\Re}$ ---called upper
(respectively lower) exhausters--- are made up of finite polytopes. In
\cite{demrya}, the authors presented a procedure ---that they called
\emph{converter}--- which permits to define from a given lower exhauster
$\underline{\Re}$ an upper exhauster $\overline{\Re}=\mathcal{F}
(\underline{\Re})$ and vice-versa (this is actually the same procedure and
coincides with the described operator $\mathcal{F}$ in the beginning of the
introduction). A lower (respectively, an upper) exhauster $\underline{\Re}$
(respectively, $\overline{\Re}$) is called \emph{stable} or \emph{reflexive},
if
\[
\underline{\Re}=\mathcal{F}\left(  \mathcal{F}(\underline{\Re})\right)
\qquad\text{(respectively, }\overline{\Re}=\mathcal{F}\left(  \mathcal{F}
(\overline{\Re})\right)  \text{).}
\]
An equivalent way to formulate the Demyanov--Ryabova conjecture is to assert
that starting with any finite (upper or lower) exhaustive family of polyhedral
functions, we eventually end up to a stable one.
\newpage

\section{Preliminary results} \label{section2}
\noindent\textbf{Notation.}

\noindent$\Re_0$ is a finite set of polytopes in $\R^N$ with $N \geq 2$. 

\noindent$S$ denotes the unit sphere of $\R^N$.

\noindent$\mathrm{ext}(\Omega)$ is the set of extreme points of a given polytope $\Omega \in \Re_0$.

\noindent$E = \bigcup_{\Omega \in \Re_0} \mathrm{ext}(\Omega)$ is the set of extreme points of all polytopes in 
$\Re_0$.

\noindent$R:=\mathrm{card}(E)$

\noindent$C:=\mathrm{conv}(E)$

\noindent We assume throughout the paper that $E$ satisfies the assumption (H1) of Theorem \ref{hyp}. For the proof of Theorem \ref{hyp}, we shall need the two following notions.
\begin{defi}
~
\begin{itemize}[leftmargin=0mm]
\item \textbf{($d$-compatible enumeration)} An enumeration $\{x_{i}\}_{i=1}^{R}$ of $E$ is called \textbf{$d$-compatible} with respect to a direction $d\in S$, provided
\begin{equation}
\langle x_{1},d\rangle\leq\langle x_{2},d\rangle\leq\cdots\leq\langle
x_{R},d\rangle. \label{d-compat}
\end{equation}
Notice that a $d$-compatible enumeration is not necessarily unique: indeed,
whenever $\langle x_{i},d\rangle=\langle x_{j},d\rangle$, for $1\leq i<j\leq
R$ the elements $x_{i}$ and $x_{j}$ can be interchanged in the above enumeration.

\item \textbf{(strict }$p$-\textbf{location)} A direction $d\in S$ is said to
\textbf{locate strictly} an element $\bar{x}\in E$ at the $p$-position (where
$p\in\{1,\ldots,R\}$), if there exists a $d$-compatible enumeration
$\{x_{i}\}_{i=1}^{R}$ of $E$ for which $x_{p}=\bar{x}$ and
\begin{equation*}
\ldots\leq\langle x_{p-1},d\rangle<\langle x_{p},d\rangle<\langle
x_{p+1},d\rangle\leq\ldots
\end{equation*}
In case $p=1$ (resp. $p=R$) the left strict inequality $\langle x_{p-1}
,d\rangle<\langle x_{p},d\rangle$ (resp. the right strict inequality $\langle
x_{p},d\rangle<\langle x_{p+1},d\rangle$) is vacuous. Notice further that
since $C$ is a polytope, assumption (H$_{1}$) yields that for every $\bar
{x}\in E$ the normal cone
\[
N_{C}(\bar{x})=\{d\in\mathbb{R}^{N}:\;\langle d,y-\bar{x}\rangle
\leq0,\;\forall y\in C\}
\]
of $C$ at $\bar{x}$ has nonempty interior (see \cite{Rock-livre,LHU} \emph{e.g.}), and every $d\in\,S\,\cap\,\mathrm{int\,}N_{C}(\bar{x})$ strictly
locates $\bar{x}$ in the $R$-position, under any $d$-compatible enumeration
$\{x_{i}\}_{i=1}^{R}$ of $E. $ 
\item \textbf{(selection)} A map $x \in E \mapsto \mathrm{e}_x \in S$ is called a selection if
$$\forall x \in E, \, \mathrm e_x \in S \,\cap\,\mathrm{int\,}N_C(x).$$ Thus, for every $x \in E$, $\mathrm{e}_x$ is a direction that strictly exposes $x$.
\end{itemize}
\end{defi}

We now begin a series of ``reordering results". The main goal is the following. Given a $d$-compatible enumeration of $E$ which locates an element $x$ at some position, say $i$, we construct a direction $d' \in S$ and a $d'$-compatible enumeration of $E$ which locates strictly $x$ to a possibly different position $p\geq i$. To construct such a $d'$, the general idea is to do small perturbations on $d$ using other well chosen directions. These perturbations need to be quantified and adequately controlled. We start with the following simple lemma.

\begin{lem}[Uniform control]\label{techlemma}
Let $d\in S$ and fix $x \in E \mapsto \mathrm{e}_x \in S\,\cap\,\mathrm{int\,}N_C(x)$ a selection. Then, there exist constants $M>0$ and $m>0$ such that, for every $x \in E$, the map $D_x$: $\R \to \R^N$ defined for every $t \in \R$ by $D_x(t) = d + t \mathrm e_x$ satisfies the following properties:
\begin{enumerate}[label=(\roman*)]
\item $D_x$ is continuous and $D_x(0) = d$. \vspace{-2mm}
\item For $t>0$ (respectively $t<0$) large enough in absolute value, any $(D_x(t)/\|D_x(t)\|)$-compatible enumeration $(x_i)_{i=1}^R$ of $E$  strictly locates $x$ at the $R$-position (resp. at the $1$-position). That is, for every $y\in E$, $y\neq x$: $\la x, D_x(t) \ra > \la y, D_x(t) \ra$ (resp. $\la x, D_x(t) \ra < \la y, D_x(t) \ra$). \vspace{-2mm}
\item \label{i1} For every $y_1,\, y_2 \in E$:  $|\la y_1 - y_2 , D_x(t) - D_x(0) \ra | \leq M |t|$. \vspace{-2mm}
\item \label{i2} For every $y \in E$, $y \neq x$: $|\la x - y , D_x(t) - D_x(0) \ra | \geq m |t|$. \vspace{-2mm}
\end{enumerate}
\end{lem}

\noindent\textbf{Proof.}
The first assertion is obvious. The second assertion is a simple consequence of the fact that $\mathrm e_x \in \mathrm{int\,} N_C(x)$ exposes $x$. 

Now let us prove \ref{i1}. We define $M = \max \{ \|y_1-y_2\| \, : \, y_1, \, y_2 \in E \} > 0$. Then, for every $y_1, \, y_2 \in E$,
\begin{equation*}
|\la y_1 - y_2 , D_x(t) - D_x(0) \ra | = |\la y_1 - y_2 , t \mathrm e_x \ra | \leq |t| \,  \|y_1 - y_2\| \, \|\mathrm e_x\| \leq M |t|.
\end{equation*}
In the same way we prove \ref{i2}. Define $m = \min \{ |\la x - y , \mathrm e_x \ra| \, : \, x,\, y \in E, \, x \neq y \} > 0$. Then for every $x,\, y \in E$ with $y \neq x$,
\begin{equation*}
|\la x - y , D_x(t) - D_x(0) \ra | = |t| \, |\la x - y , \mathrm e_x \ra| \geq m |t|.
\end{equation*}
\hfill $\square$
\bigskip

\begin{rem} \label{remconst}
Note that, whenever the selection $x \in E \mapsto \mathrm e_x \in S\,\cap\,\mathrm{int\,}N_C(x)$ is fixed, the constants $m$ and $M$ in the previous lemma hold for every function $D_x$ (and do not depend neither on $d$, nor on $x$).
\end{rem}
The next lemma will play a key role in the sequel.
\begin{lem}[Strict location in the very next position] \label{Lemma_partial location}
Let $\{x_{j}\}_{j=1}^{R}$ be a $d$-compatible enumeration of $E$ such that for some $1\leq i \leq R-1$ we have:
\begin{equation*}
\ldots \leq \langle x_{i-1},d \rangle < \langle x_{i},d \rangle \leq \langle
x_{i+1},d \rangle \leq \ldots
\end{equation*}
Then there exist a direction $d'\in S$ and a $d'$-compatible enumeration $\{y_{j}\}_{j=1}^{R}$ satisfying
\begin{equation*}
\{x_{1},\cdots,x_{i-1}\} \subset \{y_{1}\cdots,y_{i}\} 
\end{equation*}
and locating strictly $x_{i}$ at the $i+1$-position, that is,
\begin{equation*}
\left\{
\begin{array}
[c]{c}
y_{i+1}=x_{i}\\
\ldots\leq\langle y_{i},d'\rangle<\langle y_{i+1},d'\rangle<\langle
y_{i+2},d'\rangle\leq\ldots
\end{array}
\right.  
\end{equation*}
\end{lem}

\noindent\textbf{Proof.}
Throughout the proof, we fix $x \in E \mapsto \mathrm e_x \in S \,\cap\,  \mathrm{int\,}N_C(x)$ a selection and $m,M >0$ the universal constants given in Lemma \ref{techlemma} (\emph{c.f.} Remark \ref{remconst}). 

\textbf{Case 1:} $x_i$ is not strictly located in the $i$-position, that is the $d$-compatible enumeration $\{x_{j}\}_{j=1}^{R}$ verifies
\begin{equation*}
\ldots \leq \langle x_{i-1},d \rangle < \langle x_{i},d \rangle = \langle
x_{i+1},d \rangle \leq \ldots
\end{equation*}
An additional difficulty here is that they may exist more than one $y \in E$ such that $\langle x_{i},d \rangle = \langle y,d \rangle$ (that is $x_{i+1}$ may not be the unique point with this property). So let $k\in \{i-1,\ldots,R \}$ be the maximum index such that $\la x_i ,d \ra = \la x_k , d \ra$. Our strategy would be to do a small perturbation on $d$ with a good control in order to put $x_i$ at the $i$-position strictly. Of course this creates a new direction $d'$ together with a new ordering of elements in $E$ through $d'$. Then, we consider an element $y$ which is right after $x_i$ in the $d'$-ordering. Again, we do a small perturbation of $d'$ with a good control in order to reverse the order of $x_i$ and $y$. The key point is the uniform control of the employed perturbations ensuring that the element $x_i$ reaches the $i+1$-position and not a further position.

Let us write $a=\la x_i - x_{i-1}, d \ra > 0$, $c=M/m$ and let $\varepsilon > 0$ such that $a - 2c\ep >a/2>0$. Let us summarize our notations with the following picture

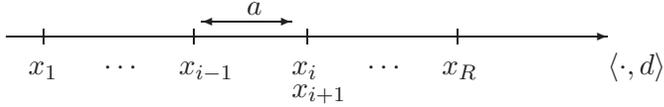
\setlength{\unitlength}{1mm}
\begin{picture}(80,14)
\put(0,5){\vector(1,0){80}}
\put(80,0){$\la \cdot, d \ra$}
\put(5,4){\line(0,1){2}} \put(3,0){$x_1$}
\put(13,0){$\cdots$}
\put(25,4){\line(0,1){2}} \put(23,0){$x_{i-1}$}
\put(40,4){\line(0,1){2}} \put(38,0){$x_{i}$} \put(38,-3){$x_{i+1}$}
\put(48,0){$\cdots$} 
\put(60,4){\line(0,1){2}} \put(58,0){$x_{R}$}
\put(26,7){\vector(1,0){12}}\put(38,7){\vector(-1,0){12}} \put(32,8){$a$}
\end{picture}

\vspace{7mm}
\noindent\textit{Step 1:} We locate $x_i$ strictly in the $i$-position but in a controlled way. Consider the map $D_{x_i}(t) = d + t \mathrm{e}_{x_i}$ defined in Lemma \ref{techlemma}, and then define the function $$\Phi: \, t \in \R \mapsto \min\limits_{j \in \{i+1,\ldots,R\} }  \la x_j - x_i, D_{x_i}(t) \ra  .$$ The map $\Phi$ is continuous, satisfies $\Phi(0)=0$ and $\lim\limits_{t \to -\infty} \Phi(t)=+\infty$. Thus, by the intermediate value theorem, there exists $t_0 <0$ such that $\Phi(t_0)= \ep$. That is $$\min\limits_{j \in \{i+1,\ldots,R\} }  \la x_j , D_{x_i}(t_0) \ra  = \la x_i , D_{x_i}(t_0) \ra + \ep.$$ 
Taking $\ep>0$ small enough we ensure that if $y \in (x_i)_{j=i+1}^R$ is such that $\la y - x_i , D_{x_i}(t_0) \ra = \ep$, then $y \in \{x_{i+1},\ldots,x_k \}$. Pick such a $y \in (x_i)_{j=i+1}^k$. Thanks to the assertion \ref{i2} of Lemma \ref{techlemma}, we have 
$$\ep = |\la y - x_i , D_{x_i}(t_0) - D_{x_i}(0) \ra| \geq m |t_0|. $$
Thus $|t_0| \leq \ep / m$. Next, thanks to the assertion \ref{i1} of Lemma \ref{techlemma}, for every $j$ in  $\{1,\ldots, i-1\}$  we have:
$$ |\la  x_i - x_j , D_{x_i}(t_0) - D_{x_i}(0) \ra | \leq M |t_0| \leq c \ep. $$
This implies that 
\begin{equation} \label{inega}
\la  x_i - x_j , D_{x_i}(t_0) \ra \geq  \la  x_i - x_j , D_{x_i}(0) \ra - c \ep \geq a - c \ep.
\end{equation}
Therefore we obtain a $(D_{x_i}(t_0)/\|D_{x_i}(t_0)\|)$-compatible enumeration $(x_i')_{i=1}^R$ satisfying $\{ x_1 , \ldots \linebreak , x_{i-1} \} = \{ x_1' , \ldots , x_{i-1}' \}$, $x_i = x_i'$ and $y=x_{i+1}' \in \{x_{i+1},\ldots,x_k\}$. We resume the situation in the following picture:

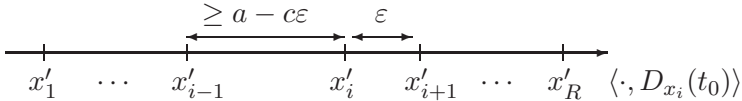
\begin{picture}(80,14)
\put(0,5){\vector(1,0){80}}
\put(80,0){$\la \cdot, D_{x_i}(t_0) \ra$}
\put(5,4){\line(0,1){2}} \put(3,0){$x_1'$}
\put(12,0){$\cdots$}
\put(24,4){\line(0,1){2}} \put(22,0){$x_{i-1}'$}
\put(45,4){\line(0,1){2}} \put(43,0){$x_{i}'$} 
\put(55,4){\line(0,1){2}}\put(53,0){$x_{i+1}'$}
\put(63,0){$\cdots$} 
\put(74,4){\line(0,1){2}} \put(72,0){$x_{R}'$}
\put(31,7){\vector(1,0){14}}\put(44,7){\vector(-1,0){20}} \put(26,9){$\geq a -c \ep$}
\put(46,7){\vector(1,0){8}}\put(54,7){\vector(-1,0){8}} \put(49,9){$\ep$}
\end{picture}

\vspace{7mm}
\noindent\textit{Step 2:} We define a new direction $\tilde{d}$ together with a $\tilde{d}$-enumeration which locates $x_i$ at the $(i+1)$-position and such that there is only one element $y$ in $E$, $y \neq x_i$, verifying $\la x_i , \tilde{d} \ra = \la y , \tilde{d} \ra$. Consider $D_{x_{i+1}'}(t) = D_{x_i}(t_0) + t \mathrm e_{x_{i+1}'}$. Reasoning as before, by the intermediate value theorem, there exists $t_1 <0$ such that $\la x_{i+1}' - x_i' , D_{x_{i+1}'}(t_1) \ra = 0$. Thanks to the assertion \ref{i2} of Lemma \ref{techlemma}, we have 
$$\ep = |\la x_{i+1}' - x_i' , D_{x_{i+1}'}(t_1) - D_{x_{i+1}'}(0) \ra| \geq m |t_1|. $$
Thus $|t_1| \leq \ep / m$. Next, thanks to the assertion \ref{i1} of Lemma \ref{techlemma}, evoking \ref{inega} under the new enumeration $\{ x_i' \}_{i=1}^R$, for every $j\in \{1,\ldots, i-1\}$ we deduce:
$$ \la  x_i' - x_j' ,D_{x_{i+1}'}(t_1) \ra  \geq \la  x_i' - x_j' , D_{x_{i+1}'}(0) \ra - M |t_1| \geq (a - c \ep ) - c \ep = a - 2c \ep. $$
Note that we also have $\la x_j' - x_{i+1}' , D_{x_{i+1}'}(t_1) \ra \geq m |t_1| $ for $j \geq i+2$. Therefore, denoting $\tilde{d} := D_{x_{i+1}'}(t_1)/\|D_{x_{i+1}'}(t_1)\|$, we may fix $(x_i'')_{i=1}^R$ a $\tilde{d}$-compatible enumeration satisfying $\{x_1'',\ldots,x_{i-1}'' \} = \{x_1,\ldots,x_{i-1} \} $, $x_i'' = x_{i+1}'$ and $ x_{i+1}'' = x_i' = x_i$. This leads us to the following configuration:

\begin{picture}(80,15)
\put(0,5){\vector(1,0){80}}
\put(80,0){$\la \cdot, D_{x_{i+1}}(t_1) \ra$}
\put(5,4){\line(0,1){2}} \put(3,0){$x_1''$}
\put(13,0){$\cdots$}
\put(25,4){\line(0,1){2}} \put(23,0){$x_{i-1}''$}
\put(45,4){\line(0,1){2}} \put(43,0){$x_{i}''$} 
\put(45,4){\line(0,1){2}}\put(43,-5){$x_{i+1}''$}
\put(60,4){\line(0,1){2}}\put(58,0){$x_{i+2}''$}
\put(66,0){$\cdots$} 
\put(74,4){\line(0,1){2}} \put(72,0){$x_{R}''$}
\put(31,7){\vector(1,0){14}}\put(44,7){\vector(-1,0){19}} \put(26,9){$\geq a -2c \ep$}
\put(46,7){\vector(1,0){13}}\put(54,7){\vector(-1,0){8}} \put(47,9){$\geq m|t_1|$}
\end{picture}

\vspace{6mm}
\noindent\textit{Step 3:} Conclusion. To complete the proof. It suffices to evoke a continuity argument and take $t_2 \in (-\infty,t_1)$ such that:
\begin{eqnarray*}
&\la& \hspace{-3mm}  x_{i}'' , D_{x_{i+1}'}(t_2) \ra > \la  x_j'' , D_{x_{i+1}'}(t_2) \ra, \, \forall \, j \in \{1,\ldots , i-1 \} \\
&\la& \hspace{-3mm} x_{i+1}'' , D_{x_{i+1}'}(t_2) \ra > \la  x_i'' , D_{x_{i+1}'}(t_2) \ra \\
&\la& \hspace{-3mm} x_{\ell}'' , D_{x_{i+1}'}(t_2) \ra >  \la  x_{i+1}'', D_{x_{i+1}'}(t_2) \ra, \, \forall \, \ell \in \{i+2,\ldots ,R \} .
\end{eqnarray*}

Setting $d' = D_{x_{i+1}'}(t_2)/\|D_{x_{i+1}'}(t_2)\|$, we deduce the existence of a $d'$-compatible enumeration $\{y_{j}\}_{j=1}^{R}$ satisfying the desired conditions. That is
$\{y_{1}\cdots,y_{i-1}\} = \{x_{1}'',\cdots,x_{i-1}''\} = \{x_{1},\cdots,x_{i-1}\}$, $y_{i}=x_{i}'' $, $y_{i+1}=x_{i+1}''=x_{i}$ and
$$\ldots\leq\langle y_{i},d'\rangle<\langle y_{i+1},d'\rangle<\langle
y_{i+2},d'\rangle\leq\ldots$$
This finishes the first part of the proof.
\bigskip

\textbf{Case 2:} $x_i$ is strictly located in the $i$-position, that is $\langle x_{i},d \rangle < \langle x_{i+1},d \rangle$. We prove that this case reduces to the first case. Indeed, consider $D_{x_i}$: $t \mapsto d +t \mathrm{e}_{x_i}$ the map given by Lemma \ref{techlemma}. Applying again the intermediate value theorem we deduce the existence of $t_0>0$ such that $$\langle x_i,D_{x_i}(t_0) \rangle = \min\limits_{j \in \{ i+1,\ldots, R \}}  \langle x_j,D_{x_i}(t_0) \rangle.$$ Thus, replacing $d$ by $\tilde{d} := D_{x_i}(t_0) / \|D_{x_i}(t_0)\| $ we obtain a $\tilde{d}$-compatible enumeration $(y_j)_{i=1}^R$ of $E$ verifying $\{x_1,\ldots,x_{i-1}\} \subset \{y_1,\ldots,y_{i-1} \}$, $y_i = x_i$ and $\langle y_{i},d' \rangle = \langle y_{i+1},d' \rangle$. Therefore we rejoined the first case. \hfill $\square$
\bigskip

\begin{rem}
It might seem strange, at a first sight, to back to the first case, since the first step of the latter was precisely to apply a perturbation that strictly locates $x_i$ in the $i$-position. However, as we pointed out in the proof, this is done in a precise quantified way.
\end{rem}
The following corollary is an easy consequence of the previous lemma and will be recalled in several occasion in the proof of Theorem \ref{hyp}.

\begin{cor}
[Reordering lemma]\label{lemorder}Let $\{x_{i}\}_{i=1}^{R}$ be a
$d$-compatible enumeration of $E$ and assume that $\la x_i , d \ra < \la x_p , d \ra$ for $1\leq i<p\leq R$. Then there exist a direction
$d^{\prime}\in S$ and a $d^{\prime}$-compatible enumeration $\{y_{j}
\}_{j=1}^{R}$ satisfying $\{x_{1},\cdots,x_{i-1}\}\subseteq\{y_{1},\cdots,y_{p-1}\}$ and strictly locating $x_i$ at the $p$-position, that is,
\begin{equation*}
\left\{
\begin{array}
[c]{c}
y_{p}=x_{i}\\
\ldots\leq\langle y_{p-1},d'\rangle<\langle y_{p},d'\rangle<\langle
y_{p+1},d'\rangle\leq\ldots
\end{array}
\right.  \label{a2}
\end{equation*}
\end{cor}

\noindent\textbf{Proof.} 
First note that if $\la x_{i-1},d \ra < \la x_i,d\ra$, then the result follows from Lemma \ref{Lemma_partial location} applied successively $p-i$ times. So let us assume that $\la x_{i-1},d \ra = \la x_i,d\ra$. Fix $\mathrm{e}_{x_i} \in \mathrm{int}N_C(x_i)$ and consider the map $D_{x_i}$: $t \in \R \mapsto d + t \mathrm{e}_{x_i}$ given by Lemma \ref{techlemma}.  Recall that $D_{x_i}$ is continuous with $D_{x_i}(0) = d$. Since $\la x_i , d \ra < \la x_p , d \ra$, there exists $t_0>0$ such that $\la x_i , D_{x_i}(t_0) \ra < \la x_p , D_{x_i}(t_0) \ra$ and $x_i$ is strictly located at some position, say $k$, in every $(D_{x_i}(t_0) /\|D_{x_i}(t_0)\|)$-compatible enumeration. Thus we set $\tilde{d} := D_{x_i}(t_0) /\|D_{x_i}(t_0)\|$ and we fix $(x_i')_{i=1}^R$ a $\tilde{d}$-compatible enumeration. Of course we have $\{x_{1},\cdots,x_{i-1}\}\subseteq\{x_{1}',\cdots,x_{k-1}'\}$ and $x_i$ is strictly located at the $k$-position in this $\tilde{d}$-compatible enumeration. Now the result follows from Lemma \ref{Lemma_partial location} applied $p-k$ times.
\hfill$\square$

\section{Proof of the main result (Theorem~\ref{hyp})}

\textbf{Extra notation}. We keep the notation introduced at the beginning of Section \ref{section2}. For the needs of the proof, we introduce some extra notation. 

\noindent Given $E_{1}\subset E$ we shall often use the abbreviate notation $[E_{1}]=\mathrm{conv}(E_{1})$. Under this notation we trivially have $C=[E]$.

\noindent Starting from a finite family of polytopes $\mathcal R_0$,  we recall that $\mathcal{R}_n = \mathcal{F}^{n}(\mathcal{R}_0)$ ($n\geq 1$) where $\mathcal{F}^{n}$ means applying the operator $\mathcal{F}$ defined in (\ref{F(R)}) $n$ times.  For $d \in  S$ and $n\geq 1$ we denote 
\[
\Omega_n(d)=[\bigcup\limits_{P\in\Re_{n-1}}\,E(P,d) \, ],
\]
where $E(P,d)=\{x\in\mathrm{ext}(P):\,\langle x,d\rangle
 =\max\langle P,d\rangle\}$. Under this notation,
\begin{equation} \label{rn}
\mathcal{R}_n  = \mathcal{F}(\mathcal{R}_{n-1}) = \{\Omega_n(d):d\in S\}.
\end{equation}
\noindent We recall that a subset $F$ of a polytope
$\Omega$ is called a face of $\Omega$ if there exists a direction $d\in S$ such that
\begin{equation}
F= \{x \in \Omega \, : \, \la x,d \ra = \min\limits_{z \in \Omega} \la z,d \ra \} . \label{face-aris}
\end{equation}
In this case we denote the face by $F(\Omega,d).$ Notice that for any $d\in S$
it holds:
\[
F(C,d)=[E(C,-d)].
\]

We are ready to proceed to the proof of Theorem~\ref{hyp}.

\bigskip

\noindent\textbf{Proof of Theorem~\ref{hyp}.} In view of Proposition
\ref{Prop-dim-1} we may assume $N\geq2$. Let us first treat the case $r_{\min
}=1,$ that is, the case where the initial family $\Re_{0}$ contains all
singletons. In this case, pick any $x\in E$ and $d\in S.$ Since $\Omega
_{x}=[x]\in\Re_{0},$ we deduce that $E(\Omega_{x},d)=\{x\}$ and consequently,
$x\in\Omega_{1}(d).$ It follows that $\Omega_{1}(d)=[E]=C$ for all $d\in S,$
that is, $\Re_{1}=\{C\}.$ Consequently, the family $\Re_{2}$ consists of all
faces of $C,$ that is,
\[
\Re_{2}=\{[E(C,d)]:d\in S\}\equiv\{F(C,d):d\in S\}.
\]
In particular, for $\bar{x}\in E$ and $d\in\mathrm{int\ }N_{C}(\bar{x})$
(direction that exposes $\bar{x}$ in $C$) we get $F(C,-d)=[\bar{x}]$, therefore
$\Re_{2}$ contains all singletons and $\Re_{3}=\{C\}=\Re_{1}$. 
\newline
Let us now treat the case $r_{\min}=R.$ In this case $\Re_{0}=\{C\}$ and we deduce, as
before, that $\Re_{1}$ is the family of all faces of $C$ and $\Re
_{2}=\{C\}=\Re_{0}$.\smallskip

It remains to treat the case $r_{\min}\notin\{1,R\}$ which is what we assume
in the sequel. In this case, we show that $\Re_{1}=\Re_{3}$ which in view of (\ref{rn}) yields $\mathcal{F}(\Re_{0})=\mathcal{F}(\Re_{2})$), \textit{i.e.}
\begin{equation}
\Omega_{1}(d)=\Omega_{3}(d)\quad\text{for every }d\in S. \label{colin}
\end{equation}
To establish (\ref{colin}) we shall proceed in three steps (Subsections
\ref{ss-1}--\ref{ss-3}), characterizing respectively, the polytopes belonging
to the families $\Re_{1}$, $\Re_{2}$ and respectively $\Re_{3}.$

\subsection{Characterization of polytopes in $\Re_{1}$.}

\label{ss-1}In this step, by means of geometric conditions on $C$ we
characterize membership of a given polytope to the family $\Re_{1}$. We start
with the biggest possible polytope, namely $C$.

\begin{prop}
\label{caracC} Assume $\Re_{0}$ satisfies (H$_{1}$), (H$_{2}$). Then the
following are equivalent:\smallskip\newline(i) $C=\Omega_{1}(d_{0})$, for some
$d_{0}\in S$ \quad(that is, $C\in\Re_{1}=\mathcal{F}(\Re_{0})$) ;\smallskip
\newline(ii) $\mathrm{card\ }\left(  F(C,d_{0})\cap E\right)  \geq r_{\min
}\quad$(that is, $C$ has a face containing at least $r_{\min}$ points).
\end{prop}

\noindent\textbf{Proof.} [(ii)$\Longrightarrow$(i)] Let us first assume that for $d_{0}\in S$ assertion
(ii) holds and let us prove that
\begin{equation}
\Omega_{1}(d_{0})=[\bigcup\limits_{\Omega\in\Re_{0}}\,E(\Omega,d_{0})]=C.
\label{a5}
\end{equation}
It suffices to prove that for each $\bar{x}\in E$ there exists a polytope
$\Omega\in\Re_{0}$ such that $\bar{x}\in
E(\Omega,d_{0}).$ Since $F(C,d_{0})$ contains at least $r_{\min}-1$ extreme
points different than $\bar{x}$, by assumption (H$_{2}$) the family $\Re_{0}$
contains the polytope $\Omega$ obtained by convexification of $\bar{x}$ and
the aforementioned $r_{\min}-1$ points. Recalling (\ref{face-aris}) we deduce
\[
E(\Omega,d_{0})=\left\{
\begin{array}
[c]{ll}
\{\bar{x}\}, & \text{if }\bar{x}\notin F(C,d_{0})\\
\Omega\cap E, & \text{if }\bar{x}\in F(C,d_{0}).
\end{array}
\right.
\]
In all cases $\bar{x}\in\,E(\Omega,d_{0})\subset\Omega_{1}(d_{0}),$ which
shows that (\ref{a5}) holds true.\smallskip

[(i)$\Longrightarrow$(ii)] Let us now assume that $C=\Omega_{1}(d_{0})$, for some $d_{0}\in S$, and let
$\{x_{i}\}_{i=1}^{R}$ be a $d_{0}$-compatible enumeration. Let $k=\max
\{i:\langle x_{i},d_{0}\rangle=\langle x_{1},d_{0}\rangle\}\ $so that
\[
F(C,d_{0})=[x_{1},\ldots,x_{k}].
\]
Assume towards a contradiction, that $k<r_{\min}$, and fix $i_{0}\in
\{1,\ldots,k\}$. Then (in view of the definition of $r_{\min}$, see (\ref{r-min})) any polytope
$\Omega\in\Re_{0}$ that contains $x_{i_{0}}$ should necessarily contain some
element $x_{j}$ with $j>k.$ In particular, $\langle x_{i_{0}},d_{0}
\rangle<\langle x_{j},d_{0}\rangle$, hence $x_{i_{0}}\notin E(\Omega,d_{0}).$
Thus $x_{i_{0}}\notin\Omega_{1}(d_{0}),$ contradicting (i).\hfill$\square$

\bigskip

Let us now characterize membership of smaller polytopes to $\Re_{1}$.

\begin{prop}
\label{carack} Assume $\Re_{0}$ satisfies (H$_{1}$), (H$_{2}$). Let $x_{1},x_{2},\ldots,x_{k}$ be distinct points in $E$ with $1 \leq k<r_{\min}$. 
The following are equivalent:\smallskip\newline(i) $C_{k}:=[E\diagdown
\{x_{1},\ldots,x_{k}\}]=\Omega_{1}(d_{k})$, for some $d_{k}\in S$ (that is,
$[E\diagdown\{x_{1},\ldots,x_{k}\}]\in\Re_{1}=\mathcal{F}(\Re_{0}
)$);\smallskip\newline(ii) There exists a $d_{k}$-compatible enumeration
$\{x_{i}^{\prime}\}_{i=1}^{R}$ of $E$ such that
\begin{equation}
\langle x_{1}^{\prime},d_k\rangle\leq\ldots\leq\langle x_{k}^{\prime}
,d_k\rangle<\langle x_{k+1}^{\prime},d_k\rangle=\cdots=\langle x_{r_{\min}
}^{\prime},d_k\rangle\leq\ldots\leq\langle x_{R},d_k\rangle, \label{aa1}
\end{equation}
and
\begin{equation}
\{x_{1},\ldots,x_{k}\}=\{x_{1}^{\prime},\ldots,x_{k}^{\prime}\}. \label{aa2}
\end{equation}
\end{prop}

\noindent\textbf{Proof.} [(ii)$\Longrightarrow$(i)] The proof is very similar
to the previous one. Let us first assume that (ii) holds for any $1\leq k<r_{\min}$ and
distinct points $x_{1},x_{2},\ldots,x_{k}\in E.$ We shall prove
\[
\bigcup\limits_{\Omega\in\Re_{0}}\,E(\Omega,d_{k})=E\diagdown\{x_{1}
,\ldots,x_{k}\},
\]
which obviously yields $\Omega_{1}(d_{k})=C_{k}.$ Pick any $i\in\{1,\ldots,k\}.$ Then by (\ref{aa2})
there exists $i_{0}\in\{1,\ldots,k\}$ with $x_{i}=x_{i_{0}}^{\prime}.$ Let
$\Omega\in\Re_{0}$ be such that $x_{i}\in\Omega.$ Then since $\mathrm{card}
\left( \Omega\right)  \geq r_{\min}>k,$ $\Omega$ should contain some
$x_{j}^{\prime}\in E$ with $\langle x_{i},d_{k}\rangle<\langle x_{j}^{\prime
},d_{k}\rangle$ (see (\ref{aa1})). Thus $x_{i}\notin E(\Omega,d_{k}).$ This
shows that
\[
\bigcup\limits_{\Omega\in\Re_{0}}\,E(\Omega,d_{k})\subset E\diagdown
\{x_{1},\ldots,x_{k}\}.
\]
Let now $\bar{x}\in E\diagdown\{x_{1},\ldots,x_{k}\}.$ Then $\langle\bar
{x},d_{k}\rangle\geq\langle x_{r_{\min}}^{\prime},d_k\rangle:=\alpha$ and by
assumption, there exist at least $\,r_{\min}-1$ extreme points with values
less or equal to $\alpha,$ forming, together with $\bar{x}$ an $r_{\min}
$-polytope $\Omega\in\Re_{0}$ for which $\bar{x}\in E(\Omega,d_{k}).$ This
shows that
\[
C_{k}:=[E\diagdown\{x_{1},\ldots,x_{k}\}]=\Omega_{1}(d_{k})\in\Re_{1},
\]
that is (i) holds.\smallskip

[(i)$\Longrightarrow$(ii)]. Assume now that for some $d_{k}\in S$ we have
$[E\diagdown\{x_{1},\ldots,x_{k}\}]=\Omega_{1}(d_{k})$, consider a $d_{k}
$-compatible enumeration $\{x_{i}^{\prime}\}_{i=1}^{R}$ of $E$, set
$\alpha:=\langle x_{r_{\min}}^{\prime},d_k\rangle$ and let $i_{1}\in
\{1,\ldots,r_{\min}\}$ (respectively, $i_{2}\in\{r_{\min},\ldots,R\}$) be the
minimum (respectively, maximum) integer $i$ such that $\langle x_{i}^{\prime
},d_k\rangle=\alpha.$ If $i_{1}=1,$ then in view of (\ref{face-aris}) the face
$F(C,d_{k})$ contains $i_{2}\geq r_{\min}$ extreme points $\{x_{1}^{\prime
},\ldots,x_{i_{2}}^{\prime}\}$. Then, according to Proposition \ref{caracC},
$\Omega_{1}(d_{k})=C$ which is a contradiction. It follows that $i_{1}>1.$ Then the
$d_{k}$-compatible enumeration satisfies
\begin{equation*}
\langle x_{1}^{\prime},d_k\rangle\leq\ldots\leq\langle x_{i_{1}-1}^{\prime
},d_k\rangle<\langle x_{i_{1}}^{\prime},d_k\rangle=\cdots=\langle x_{r_{\min}
}^{\prime},d_k\rangle\;\left( \ldots=\langle x_{i_{2}}^{\prime},d_k\rangle
< \ldots\leq\langle x_{R},d_k\rangle\right)  . 
\end{equation*}
Applying [(ii)$\Longrightarrow$(i)] for $k=i_{1}-1\in\{1,\ldots,r_{\min}-1\},$
we deduce that $C_{k}:=[E\diagdown\{x_{1}^{\prime},\ldots,x_{i_{1}-1}^{\prime
}\}],$ whence $i_{1}-1=k$ and $\{x_{1}^{\prime},\ldots,x_{i_{1}-1}^{\prime
}\}=\{x_{1},\ldots,x_{k}\}.$ The proof is complete.

\hfill$\square$

\bigskip

Let us complete this part with the following result.

\begin{prop}
\label{factminpol} Assume $\Re_{0}$ satisfies (H$_{1}$), (H$_{2}$). Then
$\Re_{1}$ does not contain any polytope of the form $[E\diagdown\{x_{1}
,\ldots,x_{k}\}]$ where $x_{1},\ldots,x_{k}\in E$ are distinct and $k\geq
r_{\min}.$
\end{prop}

\noindent\textbf{Proof.} This fact is obvious since $\Re_{0}$ contains all
possible $r_{\min}$-polytopes. In particular, there exists a polytope $\Omega$
entirely contained in $[x_{1},\ldots,x_{k}],$ and consequently for every $d\in
S$ it holds
\[
E(\Omega,d)\cap\lbrack x_{1},\ldots,x_{k}]\neq\emptyset.
\]
The proof is complete.\hfill$\square$

\bigskip

To resume the above results, we have established that a polytope $\Omega$
belongs to the family $\Re_{1}$ if and only if there is a $d_{k}$-compatible enumeration $\{x_{i}^{\prime}\}_{i=1}^{R}$ of $E$ such that 
\begin{equation*}
\Omega=[E\diagdown\{x_{1}',\ldots,x_{k}'\}]\qquad\text{(}0 \leq k<r_{\text{min}
}\text{)}
\end{equation*}
with the obvious abuse of notation: $k=0 \Longrightarrow \{x_{1}',\ldots,x_{k}'\} = \emptyset$.

\subsection{Characterization of polytopes in $\Re_{2}$.}

\label{ss-2} In this step, we shall describe the elements of the family
\[
\Re_{2}=\mathcal{F}(\Re_{1})=\{\Omega_{2}(d):d\in S\}
\]
where as usual,
\[
\Omega_{2}(d)=[\bigcup\limits_{\Omega\in\Re_{1}}\,E(\Omega,d)].
\]
Let us proceed to a complete description of the above elements. To this end,
let us fix a direction $d_{0}\in S$. By the previous step
(Subsection~\ref{ss-1}), there exists a $d_{0}$-compatible
enumeration $\{x_{i}'\}_{i=1}^{R}$ of $E$ and $k\in\{0,\ldots,r_{\min}-1\}$ such that
\begin{equation}
\Omega_{1}(d_{0})=[\bigcup\limits_{\Omega\in\Re_{0}}\,E(\Omega,d_{0}
)]=[E\diagdown\{x_{1}',\ldots,x_{k}'\}]\in\Re_{1}.\label{O1}
\end{equation}

\begin{prop}
\label{Prop-st2-1} Let $\{x_{i}'\}_{i=1}^{R}$ denote the above $d_{0}$-compatible enumeration of $E$ for which (\ref{O1}) holds. Then
\[
\Omega_{2}(-d_{0}):=[\bigcup\limits_{\Omega\in\Re_{1}}\,E(\Omega
,-d_{0})]=[x_{1}^{\prime},\ldots,x_{\ell}^{\prime}]\in\Re_{2},
\]
where
\begin{equation}
\ell=\max\{i:\langle x_{i}^{\prime},d_{0}\rangle=\langle x_{r_{\min}}^{\prime
},d_{0}\rangle\}\quad\left(\in\{r_{\min},\ldots,R\}\right).\label{O2}
\end{equation}

\end{prop}

\noindent\textbf{Proof.} Let us first assume $k\geq1.$ According to
Proposition~\ref{carack}, we have
\[\langle x_{k}^{\prime},d_{0}\rangle<\langle x_{k+1}^{\prime
},d_{0}\rangle=\langle x_{r_{\min}}^{\prime},d_{0}\rangle= \langle x_{\ell}^{\prime},d_{0}\rangle = a.
\]
Since $\Omega_{1}(d_{0})\in\Re_{1}$ the above yields
\[
E(\Omega_{1}(d_{0}),-d_{0})=\{x_{k+1}^{\prime},\ldots,x_{\ell}^{\prime}\},
\]
therefore
\[
\{x_{k+1}^{\prime},\ldots,x_{\ell}^{\prime}\}\subset\Omega_{2}(-d_{0}).
\]
Let further $m\in\{1,\ldots,k\}$ be such that
\[
F(C,d_{0})=[x_{1}^{\prime},\ldots,x_{m}^{\prime}].
\]
It follows easily that
\[
\{x_{1}^{\prime},\ldots,x_{m}^{\prime}\}=E(C,-d_{0})\subset\Omega_{2}(-d_{0}).
\]
Finally, let $i\in\{m+1,\ldots,k\}$ and let us show that $x_{i}^{\prime}
\in\Omega_{2}(-d_{0}).$ To this end, we need to exhibit a direction
$d^{\prime}\in S$ such that the polytope $\Omega_{1}(d^{\prime})\in\Re_{1}$
contains $x_{i}^{\prime}$ but does not contain any $x_{j}^{\prime}$ for $1\leq
j<i.$ (In such a case we would get $x_{i}^{\prime}\in E(\Omega_{1}(d^{\prime
}),-d_{0})\subset\Omega_{2}(-d_{0})$ and we are done.) Indeed, let $d^{\prime
}$ be given by Corollary~\ref{lemorder} for $p=r_{\min}.$ Then there exists a
$d^{\prime}$-compatible enumeration $\{y_{i}\}_{i=1}^{R}$ of $E$ locating strictly $x_i'$ in the $p=r_{\min}$ position (\textit{i.e.} $y_{r_{\min}}=x_{i}^{\prime}$) and $\{x_{1}^{\prime},\cdots,x_{i-1}^{\prime}\}\subseteq\{y_{1},\cdots,y_{r_{\min
}-1}\}$. Applying Proposition \ref{carack} [(ii) $\Longrightarrow (i)$] for $d'$ we deduce
\[
\Omega_{1}(d^{\prime})=[E\diagdown\{y_{1},\dots,y_{r_{\min}-1}\}]\in\Re_{1}
\]
and consequently
\[
x_{i}^{\prime}\in \Omega_1(d^{\prime})\quad\text{and}\quad \{x_{1}^{\prime},\cdots,x_{i-1}^{\prime}\}\cap\Omega_{1}(d^{\prime}
)=\emptyset.
\]
This proves that $\{x_{1}^{\prime},\ldots,x_{\ell}^{\prime}\}\subset\Omega_{2}
(-d_{0}).$ It remains to show that if $j>\ell$ then $x_{j}^{\prime}
\notin\Omega_{2}(-d_{0})$. Indeed, since $\ell\geq r_{\min},$ it follows from
Proposition~\ref{factminpol} that any polytope of $\Re_{1}$ should contain at
least one of the elements $\{x_{i}^{\prime}:1\leq i\leq\ell\}$. Therefore
$x_{j}^{\prime}\notin E(\Omega_{1},-d_{0})$ for all $\Omega_{1}\in\Re_{1}.$ It
follows that $\Omega_{2}(-d_{0})=[x_{1},\dots,x_{\ell}],$ as
asserted.\smallskip

Let us now assume $k=0$, that is, $\Omega_{1}(d_{0})=C.$ Then according to
Proposition~\ref{caracC} the face $F(C,d_{0})$ contains at least $r_{\min}$
points of $E$. In view of (\ref{O2}) we deduce that
\[
\lbrack x_{1}^{\prime},\ldots,x_{\ell}^{\prime}]=F(C,d_{0})=E(C,-d_{0}
)\subset\Omega_{2}(-d_{0}).
\]
Using the same argument as before,
we get that $x_{j}^{\prime}\notin\Omega_{2}(-d_{0})$ whenever $j\geq\ell +1.$ Indeed, according to Proposition~\ref{factminpol}, since $l\geq r_{min}$ any polytope of $\Re_{1}$ should contain at least one of the elements $\{x_{i}^{\prime}:1\leq i\leq\ell\}$. Thus for any polytope $\Omega_1$ in $\mathcal{R}_1$ containing $x_j$ we have $x_j \not\in E(\Omega_1,-d_0)$. The proof is complete.
\hfill$\square$
\bigskip

Since Proposition~\ref{Prop-st2-1} can be applied to all directions $d\in S$
we eventually recover a full description of polytopes in $\Re_{2}$.

\subsection{Construction of $\Re_{3}$ and conclusion.}
\label{ss-3}
In this part we prove the following assertion: For every $d\in S$, we have $\Omega_{1}(d)=\Omega_{3}(d)$. This last statement trivially implies that $\Re_{1}=\Re_{3}$ and finishes the proof of the theorem.

\noindent Let us proceed to the proof of the assertion. Fix any direction $d_{0}\in S$. According to Subsection~\ref{ss-1}, we can fix a $d_0$-compatible enumeration $(x_i')_{i=1}^R$ such that 
\[
\Omega_{1}(d_{0})=[E\diagdown\{x_{1}',\ldots,x_{k}'\}]\in\Re_{1},
\]
where $k \in \{0, \ldots, r_{\min} \}$ (under the convention that $\{x_{1}',\ldots,x_{k}' \} = \emptyset$ for $k=0$). Then, according to Proposition~\ref{carack},
\[
\Omega_{2}(-d_{0})=[x_{1}^{\prime},\dots,x_{\ell}^{\prime}]\in\Re_{2},
\]
where $\ell\geq r_{\min}$ being defined in (\ref{O2}). Thus, we are in the following configuration:
\[ \ldots \leq \la x_k', d_0 \ra  < \la x_{k+1}', d_0 \ra = \ldots =\la x_{\ell}', d_0 \ra < \la x_{\ell +1}', d_0 \ra \leq \ldots  \]

The above readily yields that
\[
E(\Omega_{2}(-d_{0}),d_{0})=\{x_{k+1}',\cdots,x_{\ell}'\}\subset\Omega_{3}
(d_{0}).
\]
Let $m\in\{\ell+1,\ldots,R\}$ be such that
\[
\lbrack x_{m}^{\prime},\ldots,x_{R}^{\prime}]=F(C,-d_{0})=E(C,d_{0}).
\]
It follows that
\[
\{x_{k+1}',\cdots,x_{\ell}'\}\cup\{x_{m}^{\prime},\cdots,x_{R}^{\prime}
\}\subset\Omega_{3}(d_{0}).
\]

Let us prove that $x_{j}'\in\Omega_{3}(d_{0})$ for all $j\in\{\ell
+1,\dots,m-1\}$. Notice that $x_j'$ is located in the $(R-j)$-position in the inverse $(-d)$-compatible enumeration. Applying Corollary \ref{lemorder} we obtain a direction $(-d')$ that pushed forward $x_j$ to the $(R - r_{\min})$-position, locating it there strictly. So we obtain a $d^{\prime}$-compatible enumeration $\{y_{j}\}_{j=1}^{R}$ such that
\[
\left\{
\begin{array}
[c]{ll}
\langle y_{R},-d^{\prime}\rangle\leq\cdots\leq\langle y_{r_{\min}},-d^{\prime
}\rangle<\langle y_{r_{\min}-1},-d^{\prime}\rangle<\cdots\leq\langle
y_{1},-d^{\prime}\rangle & \\
y_{r_{\min}}=x_{j}' & \\
\{x_{j+1}',\dots,x_{R}'\}\subseteq\{y_{r_{\min}+1},\dots,y_{R}\} &
\end{array}
\right.
\]
Writing the above assertion in reverse order yields
$$\langle y_{1},d^{\prime}\rangle  \leq  \cdots <  \langle y_{r_{\min-1}},d^{\prime}\rangle   <   \langle y_{r_{\min}},d^{\prime}\rangle < \cdots \leq \langle
y_{R},d^{\prime}\rangle .$$
It follows by Proposition~\ref{carack} [(ii) $\Longrightarrow$ (i)] that
\[
\Omega_{1}(d^{\prime})=[E\backslash\{y_{1},\dots,y_{r_{\min}-1}\}]\in\Re_{1}
\]
and consequently, $y_{r_{\min}} = x_{j}'\in E(\Omega_{1}(d^{\prime}),d_{0})\subset\Omega
_{3}(d_{0}).$ 

It remains to prove that $x_{j}'\not \in \Omega_{3}(d_{0})$
whenever $j\in\{1,\dots,k\}$. Indeed, if this were not the case, then there
would exist a polytope $\Omega\in\Re_{2}$ such that $x_{j}'\in E(\Omega,d_{0})$
and consequently the polytope $\Omega$ cannot contain any other element $x\in
E$ with $\langle x,d_{0}\rangle>\langle x_{j}',d_{0}\rangle$. In particular
$\{x_{k+1}^{\prime},\dots,x_{R}'\}\cap\Omega=\emptyset$. Thus such a polytope could contain at most $k$ points of $E$ with $k < r_{min}$, which is impossible
according to Proposition~\ref{Prop-st2-1} (every polytope of $\Re_2$ contains at least $r_{min}$ points of $E$). It follows that
\[
\Omega_{3}(d_{0})=[E\backslash\{x_{1}',\dots,x_{k}'\}]=\Omega_{1}(d_{0}),
\]
which proves the assertion and the theorem.
\hfill$\square$
\bigskip

\subsection{Weakening assumption (H2)}

A careful inspection of the previous proof reveals that some $r_{\min}
$-polytopes do not intervene in the construction of the family $\Re
_{1}=\emph{F}(\Re_{0})$ and consequently assumption $(H_{2})$ can be relaxed
as follows (we leave the details to the reader):

\begin{itemize}
\item[$\mathrm{(H}_{2}^{\prime}\mathrm{)}$] The family $\Re_{0}$ contains all
$r_{\min}$-polytopes of the form $[x_{1},\ldots,x_{r_{\min}}]$ for which there exists a
direction $d\in S$ and a $d$-compatible enumeration $\{x_{i}^{\prime}
\}_{i=1}^{R}$ such that
\[
\{x_{1},\ldots,x_{r_{\min}}\}=\{x_{1}^{\prime},\ldots,x_{r_{\min}}^{\prime
}\}\text{\quad and\quad}\langle x_{r_{\min}}^{\prime},d\rangle<\langle
x_{r_{\min}+1}^{\prime},d\rangle.
\]

\end{itemize}

\bigskip

\textbf{Acknowledgments.} The authors thank Vera Roshchina (University
of Melbourne) for introducing to them the conjecture. They also thank Abderrahim
Hantoute (CMM, University of Chile) and Bernard Baillon (University Paris~1) for useful discussions.
Major part of this work has been accomplished during a research visit of the second author to the
Department of Mathematical Engineering of the University of Chile (October
2016) and of the first author to the Laboratory of Mathematics of the
University of Franche-Comt\'{e} in Besan\c{c}on (December 2016). The
authors thank their hosts for hospitality.

\vspace{1cm}

\noindent Aris Daniilidis

\smallskip

\noindent DIM--CMM, UMI CNRS 2807\newline Beauchef 851 (Torre Norte, piso~5),
Universidad de Chile, $\ldots$ , Santiago, Chile \smallskip

\noindent E-mail: \texttt{arisd@dim.uchile.cl} \newline\noindent
\texttt{http://www.dim.uchile.cl/{\raise.17ex\hbox{$\scriptstyle\sim$}}arisd}

\smallskip

\noindent Research supported by the grants ECOS C14E06, BASAL PFB-03, FONDECYT
1171854 (Chile) and MTM2014-59179-C2-1-P (MINECO-ERDF, Spain and EU).

\bigskip

\noindent Colin Petitjean

\smallskip

\noindent Universit\'e Franche-Comt\'e, Laboratoire de Math\'ematiques UMR
6623\newline16 route de Gray, F-25030 Besan\c con Cedex, France

\noindent E-mail: \texttt{colin.petitjean@univ-fcomte.fr} \newline\noindent
\texttt{http://cpetit13.perso.math.cnrs.fr/}

\smallskip

\noindent Research supported by the grant ECOS-Sud C14E06 (France).


\begin{thebibliography}{9}



\bibitem {demrya}\textsc{V. F. Demyanov, J. A. Ryabova,} \textit{Exhausters,
coexhausters and converters in nonsmooth analysis}, \emph{Discrete Contin.
Dyn. Syst.} \textbf{31} (2011), no. 4, 1273--1292

\bibitem {rubinov}\textsc{V. F. Demyanov et al.} \textit{Negladkie zadachi
teorii optimizatsii i upravleniya}, (Leningrad. Univ. 1982).

\bibitem {LHU}\textsc{J.-B. Hiriart-Urruty, C. Lemar\'{e}chal,}
\emph{Fundamentals of convex analysis}, Grundlehren Text Editions. (Springer, 2001).

\bibitem {phelps}\textsc{R. R. Phelps}, \textit{Convex functions, monotone
operators and differentiability}, Second, Lectures Notes in Mathematics,
\textbf{vol. 1364}, (Springer, 1993).

\bibitem {psh}\textsc{B. N. Pshenichnyi,} \textit{Convex Analysis and Extremal
Problems}, (Nauka, Moscou, 1980).

\bibitem {Rock-livre}\textsc{T. R. Rockafellar,} \emph{Convex analysis},
Princeton Mathematical Series, No. 28, (Princeton, N.J. 1970).

\bibitem {sang}\textsc{T. Sang,} \textit{On the conjecture by Demyanov-Ryabova
in converting finite exhausters}, ArXiv e-prints (2016).
\end{thebibliography}
\end{document}